\documentclass[12pt, twoside, leqno]{article}
\usepackage{amsmath,amsthm}
\usepackage{amssymb,latexsym}
\usepackage{enumerate}

\newtheorem{thm}{Theorem}[section]
\newtheorem{prop}[thm]{Proposition}
\newtheorem{cor}[thm]{Corollary}
\newtheorem{lem}[thm]{Lemma}

\theoremstyle{definition}

\numberwithin{equation}{section}

\usepackage{mathrsfs}
\usepackage{CJK}
\usepackage{amsmath}
\usepackage{fancyhdr}
\usepackage{enumerate}

\begin{document}

\baselineskip=17pt

\title{Some polynomial maps with Jacobian rank two or three}
\author{ Dan Yan \footnote{ The author is supported by the National Natural
Science Foundation of China (Grant No.11601146; 11871241), the
Natural Science Foundation of Hunan Province (Grant No.2016JJ3085)
and the Construct
Program of the Key Discipline in Hunan Province.}\\
MOE-LCSM,\\ School of Mathematics and Statistics,\\
 Hunan Normal University, Changsha 410081, China \\
\emph{E-mail:} yan-dan-hi@163.com \\
}
\date{}

\maketitle

\renewcommand{\thefootnote}{}

\renewcommand{\thefootnote}{\arabic{footnote}}
\setcounter{footnote}{0}

\begin{abstract} In the paper, we first classify all
polynomial maps of the form $H=(u(x,y,z),v(x,y,z), h(x,y))$ in the
case that $JH$ is nilpotent and $\deg_zv\leq 1$.  After that, we
generalize the structure of $H$ to
$H=\big(H_1(x_1,x_2,\ldots,x_n),\allowbreak
b_3x_3+\cdots+b_nx_n+H_2^{(0)}(x_2),H_3(x_1,x_2),\ldots,H_n(x_1,x_2)\big)$.
\end{abstract}
 \vskip 2.5mm

\section{Introduction}

Throughout this paper, we will write $K$ for any field with
characteristic zero and $K[x]=K[x_1,x_2,\ldots,x_n]$ for the
polynomial algebra over $K$ with $n$ indeterminates. Let
$F=(F_1,F_2,\ldots,F_n):K^n\rightarrow K^n$ be a polynomial map,
that is, $F_i\in K[x]$ for all $1\leq i\leq n$. Let
$JF=(\frac{\partial F_i}{\partial x_j})_{n\times n}$ be the Jacobian
matrix of $F$. For $H_i, u_{i-1}\in K[x]$, we abbreviate $\frac{\partial
H_i}{\partial x_j}$ as $H_{ix_j}$ and $\frac{\partial
u_{i-1}}{\partial x_j}$ as $u_{(i-1)x_j}$, and define $\deg_{x_i} f$ as the
highest degree of variable $x_i$ in $f$. $P_n(i,j)$ denotes the $n
\times n$ elementary permutation matrix which interchanges
coordinates $i$ and $j$, and $P_n(i(a),j)$ denotes the $n
\times n$ elementary matrix which add $a$ times of the $i$-th row to the $j$-th row.

The Jacobian Conjecture (JC) raised by O.H. Keller in 1939 in
\cite{1} states that a polynomial map
$F: K^n\rightarrow K^n$ is invertible if the Jacobian
determinant $\det JF$ is a nonzero constant. This conjecture has
been attacked by many people from various research fields, but it is
still open, even for $n\geq 2$. Only the case $n=1$ is obvious. For
more information about the wonderful 70-year history, see \cite{2},
\cite{3}, and the references therein.

In 1980, S.S.S.Wang (\cite{4}) showed that the JC holds for all
polynomial maps of degree 2 in all dimensions (up to an affine
transformation). A powerful result is the reduction to degree
3, due to H.Bass, E.Connell and D.Wright (\cite{2}) in 1982 and
A.Yagzhev (\cite{5}) in 1980, which asserts that the JC is true if
the JC holds for all polynomial maps $x+H$, where $H$ is homogeneous
of degree 3. Thus, many authors study these maps and led to pose the
following problem.

 {\em (Homogeneous) dependence problem.} Let $H=(H_1,\ldots,H_n)\in
K[x]$ be a (homogeneous) polynomial map of degree $d$ such
that $JH$ is nilpotent and $H(0)=0$. Whether $H_1,\ldots,H_n$ are
linearly dependent over $K$?

The answer to the above problem is affirmative if rank$JH\leq 1$
(\cite{2}). In particular, this implies that the Dependence Problem
has an affirmative answer in the case $n=2$. M. de Bondt and Van den
Essen give an affirmative answer to the above problem in the case that $H$
is homogeneous and $n=3$ (\cite{8}).

With restrictions on the degree of $H$, more positive results are known.
For cubic homogeneous $H$, the case $n = 4$ has been solved affirmatively
by Hubbers in \cite{7}, using techniques of \cite{6}. For cubic homogeneous
$H$ with rank$JH = 2$, the Dependence Problem has an affirmative answer for every
$n$, because the missing case $n \geq 5$ follows from \cite[Theorem 4.3.1]{HKM}.
For cubic $H$, the case $n = 3$ has been solved affirmatively as well,
see e.g.\@ \cite[Corollary 4.6.6]{HKM}).

For quadratic $H$, the Dependence Problem has an affirmative answer if
rank$\allowbreak JH \leq 2$ (see \cite{B2} or \cite[Theorem 3.4]{12}), in
particular if $n \leq 3$. For quadratic homogeneous $H$,
the Dependence Problem has an affirmative answer in the case $n \leq 5$,
and several authors contributed to that result. See \cite[Appendix A]{HKM}
and \cite{XS5} for the case $n = 5$.

The first counterexamples to the Dependence Problem were found by
Van den Essen (\cite{9}, \cite[Theorem 7.1.7 (ii)]{3}). He
constructs counterexamples for all $n \geq 3$. In another paper
(\cite{E}), he constructs a quadratic counterexample for $n = 4$, which can be
generalized to arbitrary even degree (see \cite[Example 8.4.4]{3} for degree $4$.

 M. de Bondt was the
first who found homogeneous counterexamples (\cite{10}). He
constructed homogeneous counterexamples of $6$
for $n = 5$, homogeneous counterexamples of degree $4$ and $5$
for all $n \geq 6$, and cubic homogeneous counterexamples for all $n \geq 10$.
Homogeneous counterexamples of larger degrees can be made as well,
except for $n = 5$ and odd degrees. A cubic homogeneous counterexample
for $n = 9$ can be found in \cite{SFGZ}, see also \cite[Section 4.2]{HKM}.

In \cite{18}, Chamberland and Van den Essen classified
all polynomial maps of the form
$$H=\big(u(x_1,x_2),v(x_1,x_2,x_3),h(u(x_1,x_2),v(x_1,x_2,x_3))\big)$$
with $JH$ nilpotent. The author and Tang \cite{13} classified
all polynomial maps of the form
$H=\big(u(x_1,x_2),v(x_1,x_2,x_3),h(x_1,x_2,x_3)\big)$ with some
conditions. In \cite{14}, the author and M. de Bondt classify all polynomial maps
of the form
$$H=\big(H_1(x_1,x_2,\ldots,x_n),H_2(x_1,x_2),H_3(x_1,x_2,H_1),\ldots,H_n(x_1,x_2,H_1)\big)$$
with $JH$ nilpotent.
Casta\~{n}eda and Van den Essen classify in \cite{CE} all polynomial maps
of the form
$$H=\big(u(x_1,x_2),u_2(x_1,x_2,x_3),u_3(x_1,x_2,x_4),\ldots,u_{n-1}(x_1,x_2,x_n),
h(x_1,x_2)\big)$$ with $JH$ nilpotent.

A polynomial map of the form
$(x_1,\ldots,x_{i-1},x_i+Q,x_{i+1},\ldots,x_n)$ is \emph{elementary} if
$Q\in K[x_1,\ldots,x_{i-1},x_{i+1},\ldots,x_n]$. A polynomial map is
called \emph{tame} if it is a finite composition of invertible linear maps
and elementary maps.

In the paper, we first classify all polynomial maps of the form
$H=(u(x,y,z),\allowbreak v(x,y,z), h(x,y))$ in the case that $JH$ is nilpotent and $\deg_zv\leq 1$. Then, in section
3, we extend these results to the case where
$$H=\big(H_1(x_1,x_2,\ldots,x_n),\allowbreak
b_3x_3+\cdots+b_nx_n+H_2^{(0)}(x_2),H_3(x_1,x_2),\ldots,H_n(x_1,x_2)\big).$$
In particular, we prove that $F=x+H$ is tame.

\section{Polynomial maps of the form $H=(u(x,y,z),v(x,\allowbreak y,z), h(x,y))$}

In the section, we classify polynomial maps of the form
$H=(u(x,y,z),v(x,y,z),\allowbreak h(x,y))$ in the case that $JH$ is
nilpotent, $\deg_zv(x,y,z)\leq 1$. Firstly, we prove in Lemma 2.2
that $u,v,h$ are linearly dependent in the case that $JH$ is
nilpotent, $\deg_zv(x,y,z)=1$ and $\deg_zu=2$. Then we prove
that $u,v,h$ are linearly dependent in Theorem 2.3 in the case that
$JH$ is nilpotent and $\deg_zv=1$ and $\deg_zu> 2$.

\begin{lem} \label{lem2.1}
Let $u,v\in K[x,y,z]$, $u=u_dz^d+\cdots+u_1z+u_0$,
$v=v_lz^l+\cdots+v_1z+v_0$ with $u_dv_l\neq 0$. If $u_x+v_y=0$ and
$l\leq d$, then we have the following equations
\begin{equation}\label{eq2.1}
u_{dx}=\cdots=u_{(l+1)x}=0
\end{equation}

 and
\begin{equation}\label{eq2.2}
u_{ix}+v_{iy}=0
\end{equation}

for $0\leq i\leq l$.
\end{lem}
\begin{proof}
We have the conclusion by comparing the coefficients of $z^j$ of the
equation $u_x+v_y=0$ for $0\leq j\leq d$.
\end{proof}

\begin{lem}\label{lem2.2}
Let $H=(u(x,y,z),v(x,y,z),h(x,y))$ be a polynomial map with
$\deg_zv(x,y,z)=1$. Assume that $H(0)=0$ and $\deg_zu=2$. If $JH$ is
nilpotent, then $u,v,h$ are linearly dependent.
\end{lem}
\begin{proof}
Since $JH$ is nilpotent, we have the following equations:
\begin{eqnarray}
  u_x+v_y = 0,\label{eq2.3}\\
  u_xv_y-v_xu_y-h_xu_z-h_yv_z=0,\label{eq2.4}\\
  v_xh_yu_z-h_xv_yu_z+h_xu_yv_z-u_xh_yv_z = 0.\label{eq2.5}
  \end{eqnarray}
Let $u,v$ be as in Lemma \ref{lem2.1}. Since $\deg_zu=2$, $\deg_zv=1$, it
follows from equation \eqref{eq2.3} and Lemma \ref{lem2.1} that
\begin{equation}\label{eq2.6}
u_{2x}=0
\end{equation}
and
\begin{equation}\label{eq2.7}
u_{ix}+v_{iy}=0
\end{equation}
for $0\leq i\leq 1$. It follows from equations \eqref{eq2.4} and \eqref{eq2.6} that
\begin{equation}\label{eq2.8}
\begin{split}
(u_{1x}z+u_{0x})(v_{1y}z+v_{0y})-(v_{1x}z+v_{0x})(u_{2y}z^2+u_{1y}z+u_{0y})\\
-h_x(2u_2z+u_1)-h_yv_1=0.
\end{split}
\end{equation}
We always view that the polynomials are in $K[x,y,z]$ with coefficients in $K[x,y]$ while we compare the coefficients of the degree of $z$.
Comparing the coefficients of $z^3$ and $z^2$ of the above equation,
 we have $v_{1x}u_{2y}=0$ and
 \begin{equation}\label{eq2.9}
u_{1x}v_{1y}-v_{1x}u_{1y}-v_{0x}u_{2y}=0.
\end{equation}
Thus, we have $u_{2y}=0$ or $v_{1x}=0$.\\

(I) If $u_{2y}=0$, then we have $u_2\in K^*$ and
\begin{equation}\label{eq2.10}
u_{1x}v_{1y}-v_{1x}u_{1y}=0
\end{equation}
by equations \eqref{eq2.6} and \eqref{eq2.9} respectively.
It follows from equation \eqref{eq2.7} that $u_{1x}=-v_{1y}$. Thus, there
exists $P\in K[x,y]$ such that $u_1=P_y$, $v_1=-P_x$. It
follows from equation \eqref{eq2.10}) that $P_{xy}^2-P_{xx}P_{yy}=0$. Then it
follows from Lemma 2.1 in \cite{18} that
\begin{equation}\label{eq2.11}
u_1=P_y=bf(ax+by)+c_2
\end{equation}
and
\begin{equation}\label{eq2.12}
v_1=-P_x=-af(ax+by)+c_1
\end{equation}
for some $f(t)\in K[t]$ and $f(0)=0$, $a,b\in K^*$,
$c_1,c_2\in K$. Then we have the following equations:
\begin{equation}\label{eq2.13}
u_{1x}v_{0y}+u_{0x}v_{1y}-v_{1x}u_{0y}-v_{0x}u_{1y}-2u_2h_x=0
\end{equation}
and
\begin{equation}\label{eq2.14}
u_{0x}v_{0y}-v_{0x}u_{0y}-u_1h_x-v_1h_y=0
\end{equation}
by comparing the coefficients of $z$ and $z^0$ of equation \eqref{eq2.8} respectively.
It follows from equations \eqref{eq2.5} and \eqref{eq2.6} that
\begin{equation}\label{eq2.15}
\begin{split}
[(v_{1x}z+v_{0x})h_y-h_x(v_{1y}z+v_{0y})](2u_2z+u_1)\\
+[h_x(u_{1y}z+u_{0y})-h_y(u_{1x}z+u_{0x})]v_1=0.
\end{split}
\end{equation}
Comparing the coefficients of $z^2,z,z^0$ of equation \eqref{eq2.15}, we
have the following equations:
\begin{equation}\label{eq2.16}
v_{1x}h_y-h_xv_{1y}=0,
\end{equation}
\begin{equation}\label{eq2.17}
(v_{0x}h_y-h_xv_{0y})2u_2+(v_{1x}h_y-h_xv_{1y})u_1+(h_xu_{1y}-h_yu_{1x})v_1=0
\end{equation}
and
\begin{equation}\label{eq2.18}
(v_{0x}h_y-h_xv_{0y})u_1+(h_xu_{0y}-h_yu_{0x})v_1=0.
\end{equation}
It follows from equations \eqref{eq2.12} and \eqref{eq2.16} that $af'\cdot
(bh_x-ah_y)=0$. Thus, we have $a=0$ or $f'=0$ or $bh_x=ah_y$.\\

(i) If $a=0$, then
\begin{equation}\label{eq2.19}
u_1=bf(by)+c_2,~~~v_1=c_1\in  K^*
\end{equation}
It follows from equations \eqref{eq2.13} and \eqref{eq2.19} that $2u_2h_x=-b^2f'(by)v_{0x}$.
Integrating with respect to $x$ of two sides of the above equation,
we have
\begin{equation}\label{eq2.20}
h=-\frac{b^2}{2u_2}f'(by)v_0+\frac{c(y)}{2u_2}
\end{equation}
for some $c(y)\in K[y]$. Substituting equations \eqref{eq2.19} and \eqref{eq2.20}
 into equation \eqref{eq2.17}, we have the following equation:
$$v_{0x}[-b^3f''(by)v_0+c'(y)-\frac{b^4}{2u_2}v_1\cdot (f'(by))^2]=0.$$
Thus, we have $v_{0x}=0$ or $f''(by)=0$ and
$c'(y)=\frac{b^4}{2u_2}v_1\cdot (f'(by))^2$.

If $v_{0x}=0$, then it follows from equations \eqref{eq2.13} and \eqref{eq2.19} that
$h_x=0$. It follows from equation \eqref{eq2.18} that $v_1h_yu_{0x}=0$.
Thus, we have $h_y=0$ or $u_{0x}=0$. If $u_{0x}=0$, then it follows
from equation \eqref{eq2.14} that $h_y=0$. If $h_y=0$, then $h=0$ because
$h(0,0)=0$. Thus, $u,v,h$ are linearly dependent.

If $f''(by)=0$, then $f'(by)\in K$ and
$c'(y)=\frac{b^4}{2u_2}v_1(f'(by))^2\in K$. Let $l:=b^2f'(by)$
and $c:=\frac{l^2}{(2u_2)^2}v_1$. Then it follows from equation \eqref{eq2.19} that $v_1\in  K^*$ and $u_1=ly+c_2$ for some $c_2\in K$. Since $h(0,0)=0$, it follows from equation \eqref{eq2.20} that
\begin{equation}\label{eq2.21}
h=-\frac{l}{2u_2}v_0+c\cdot y.
\end{equation}
It follows from equations \eqref{eq2.14} and \eqref{eq2.21} that
\begin{equation}\label{eq2.22}
v_{0x}u_{0y}-u_{0x}v_{0y}=\frac{l}{2u_2}u_1v_{0x}+\frac{l}{2u_2}v_1v_{0y}-cv_1.
\end{equation}
It follows from equations \eqref{eq2.18} and \eqref{eq2.21} that
\begin{equation}\label{eq2.23}
c\cdot u_1v_{0x}+v_1[-\frac{l}{2u_2}(v_{0x}u_{0y}-u_{0x}v_{0y})-c\cdot u_{0x}]=0.
\end{equation}
Substituting equation \eqref{eq2.22} into equation \eqref{eq2.23}, we have the
following equation:

$$c\cdot u_1v_{0x}+v_1[-\frac{l^2}{(2u_2)^2}u_1v_{0x}+\frac{lc}{2u_2}v_1-\frac{l^2}{(2u_2)^2}v_1v_{0y}-c\cdot u_{0x}]=0.$$

Since $c=\frac{l^2}{(2u_2)^2}v_1$, the above equation has the
following form:
$$\frac{lc}{2u_2}v_1-c(v_{0y}+u_{0x})=0.$$
Substituting equation \eqref{eq2.7}$(i=0)$ into the above equation, we have
$\frac{lc}{2u_2}v_1=0$. That is, $lc=0$. Since
$c=\frac{l^2}{(2u_2)^2}v_1$, we have $c=l=0$. It follows from equation
\eqref{eq2.21} that $h=0$. Thus, $u,v,h$ are linearly dependent.\\

(ii) If $f'=0$, then $f=0$ because $f(0)=0$. That is, $u_1=c_2$,
$v_1=c_1\in  K^*$. It follows from equation \eqref{eq2.13} that
$h_x=0$. It follows from equation \eqref{eq2.17} that $v_{0x}h_y=0$. Thus,
we have $h_y=0$ or $v_{0x}=0$.\\
If $h_y=0$, then $h=0$ because $h(0,0)=0$. Thus, $u,v,h$ are
linearly dependent.\\
If $v_{0x}=0$, then it follows from equation \eqref{eq2.18} that
$u_{0x}h_y=0$. That is, $u_{0x}=0$ or $h_y=0$. If $u_{0x}=0$, then
it follows from equation \eqref{eq2.14} that $h_y=0$. Thus, we have that
$h_y=0$. It reduces to the above case.\\

(iii) If $bh_x=ah_y$, then we can assume that $a\cdot f'\neq 0$, hence we have
\begin{equation}\label{eq2.24}
h_y=\frac{b}{a}h_x.
\end{equation}
Let $\bar{x}=ax+by$, $\bar{y}=y$. It follows from equation \eqref{eq2.24} that we
have $h_{\bar{y}}=0$. That is, $h\in  K[ax+by]$. It
follows from equations \eqref{eq2.11}, \eqref{eq2.12}, \eqref{eq2.17} and \eqref{eq2.24} that
\begin{equation}\label{eq2.25}
h_x(\frac{b}{a}v_{0x}-v_{0y})=0.
\end{equation}
It follows from equations \eqref{eq2.18} and \eqref{eq2.24} that
\begin{equation}\label{eq2.26}
u_1h_x(\frac{b}{a}v_{0x}-v_{0y})+v_1h_x(u_{0y}-\frac{b}{a}u_{0x})=0.
\end{equation}
It follows from equations \eqref{eq2.25} and \eqref{eq2.26} that $h_x=0$ or
$bv_{0x}=av_{0y}$ and $bu_{0x}=au_{0y}$. \\
If $h_x=0$, then it follows from equation \eqref{eq2.24} that $h_y=0$. Thus,
we have $h=0$ because $h(0,0)=0$. so $u,v,h$ are linearly
dependent.\\
If $bv_{0x}=av_{0y}$ and $bu_{0x}=au_{0y}$, then $v_0,u_0\in
K[ax+by]$ by the same reason as $h$. Thus, it follows from equations \eqref{eq2.11},  \eqref{eq2.12} and
 \eqref{eq2.13} that $h_x=0$. It reduces to the former case.\\

(II) If $v_{1x}=0$, then it follows from equation \eqref{eq2.9} that
\begin{equation}\label{eq2.27}
u_{1x}v_{1y}-v_{0x}u_{2y}=0.
\end{equation}
It follows from equation \eqref{eq2.8} that
\begin{equation}\label{eq2.28}
(u_{1x}z+u_{0x})(v_{1y}z+v_{0y})-v_{0x}(u_{2y}z^2+u_{1y}z+u_{0y})-h_x(2u_2z+u_1)-h_yv_1=0.
\end{equation}
Comparing the coefficients of $z^2,~z,~z^0$ of equation \eqref{eq2.28}, we
have the following equations:
\begin{equation}\label{eq2.29}
u_{1x}v_{1y}-v_{0x}u_{2y} = 0,
\end{equation}
\begin{equation}\label{eq2.30}
u_{1x}v_{0y}+u_{0x}v_{1y}-v_{0x}u_{1y}-2u_2h_x = 0,
\end{equation}
\begin{equation}\label{eq2.31}
 u_{0x}v_{0y}-v_{0x}u_{0y}-u_1h_x-v_1h_y = 0.
 \end{equation}
It follows from equations \eqref{eq2.5} and \eqref{eq2.6} that
$$[v_{0x}h_y-(v_{1y}z+v_{0y})h_x](2u_2z+u_1)+[h_x(u_{2y}z^2+u_{1y}z+u_{0y})-h_y(u_{1x}z+u_{0x})]v_1=0.$$
Comparing the coefficients of $z^2,~z$ and $z^0$ of the above
equation, we have the following equations:
\begin{equation}\label{eq2.32}
h_x(v_1u_{2y}-2u_2v_{1y}) = 0,
\end{equation}
\begin{equation}\label{eq2.33}
-v_{1y}h_xu_1+2u_2(v_{0x}h_y-v_{0y}h_x)+v_1(h_xu_{1y}-h_yu_{1x}) = 0
\end{equation}
and
\begin{equation}\label{eq2.34}
u_1(v_{0x}h_y-v_{0y}h_x)+v_1(h_xu_{0y}-h_yu_{0x}) = 0.
\end{equation}
It follows from equation \eqref{eq2.32} that $h_x=0$ or
$v_1u_{2y}=2u_2v_{1y}$.

If $h_x=0$, then it follows from equation \eqref{eq2.33} that
$h_y(2u_2v_{0x}-u_{1x}v_1)=0$. Thus, we have that $h_y=0$ or
$2u_2v_{0x}=v_1u_{1x}$. If $h_y=0$, then $h=0$ because $h(0,0)=0$.
Thus, $u,v,h$ are linearly dependent. If $2u_2v_{0x}=v_1u_{1x}$,
then it follows from equation \eqref{eq2.7}(i=1) that
\begin{equation}\label{eq2.35}
2u_2v_{0x}=-v_1v_{1y}.
\end{equation}
Substituting equations \eqref{eq2.35} and \eqref{eq2.7} to equation \eqref{eq2.29} , we have the following equation: $v_{1y}(2u_2v_{1y}-v_1u_{2y})=0$. Thus, we
have $v_{1y}=0$ or $2u_2v_{1y}=v_1u_{2y}$. \\
If $v_{1y}=0$, then it follows from equation \eqref{eq2.29} that
$v_{0x}u_{2y}=0$. Thus, we have $v_{0x}=0$ or $u_{2y}=0$. If
$u_{2y}=0$, then it reduces to (I). If $v_{0x}=0$, then it
follows form equation \eqref{eq2.34} that $h_yu_{0x}=0$. Thus, we have
$h_y=0$ or $u_{0x}=0$. If $u_{0x}=0$, then it follows from equation \eqref{eq2.31} that $h_y=0$. Therefore, we have $h=0$ because $h(0,0)=0$. Thus, $u,v,h$ are linearly dependent.

If $2u_2v_{1y}=v_1u_{2y}$, then we have
\begin{equation}\label{eq2.36}
\frac{u_{2y}}{u_2}=2\frac{v_{1y}}{v_1}.
\end{equation}
Suppose that $u_{2y}v_{1y}\neq 0$. Then we have
$$u_2=e^{\bar{c}(x)}v_1^2$$
by integrating the two sides of \eqref{eq2.36} with respect to $y$.
where $\bar{c}(x)$ is a function of $x$. Since $u_2,v_1\in
K[x,y]$, we have $e^{\bar{c}(x)}\in  K(x)$. That is,
$u_2=c(x)v_1^2$, where $c(x)$ is not equal to zero and belongs to $K(x)$. Let
$c(x)=\frac{c_1(x)}{c_2(x)}$ with $c_1(x), c_2(x)\in K[x]$ and
$c_1(x)\cdot c_2(x)\neq 0$. Then it follows from equations \eqref{eq2.29}
and \eqref{eq2.7} that
$$v_{1y}(2c(x)v_1v_{0x}+v_{1y})=0.$$
That is,
\begin{equation}\label{eq2.37}
2c_1(x)v_1v_{0x}=-c_2(x)v_{1y}.
\end{equation}
If $v_{0x}\neq 0$, then we have that $v_{1y}=0$ by comparing the
degree of $y$ of equation \eqref{eq2.37}. Thus, we have $v_{0x}=0$. It follows from equation \eqref{eq2.37} that $v_{1y}=0$. This contradicts with our assumption.
Therefore, we have $v_{1y}=v_{0x}=0$. It follows
from equation \eqref{eq2.36} that
$u_{2y}=v_{1y}=0$. which reduces to (I).
\end{proof}

\begin{thm}\label{thm2.3}
Let $H=(u(x,y,z),v(x,y,z),h(x,y))$ be a polynomial map with
$\deg_zv(x,y,z)=1$. Assume that $H(0)=0$ and $\deg_zu\geq 2$. If
$JH$ is nilpotent, then $u,v,h$ are linearly dependent.
\end{thm}
\begin{proof}
Let $u$, $v$ be as in Lemma \ref{lem2.1}.

If $\deg_zu=2$, then the conclusion follows from Lemma \ref{lem2.2}.

If $\deg_zu\geq 3$, then it follows from equation \eqref{eq2.3} and Lemma \ref{lem2.1} that
\begin{equation}\label{eq2.38}
u_{dx}=\cdots=u_{2x}=0
\end{equation}
and equation \eqref{eq2.7} is true. It follows from equations \eqref{eq2.4} and \eqref{eq2.38} that
\begin{equation}\label{eq2.39}
\begin{split}
(u_{1x}z+u_{0x})(v_{1y}z+v_{0y})-(v_{1x}z+v_{0x})(u_{dy}z^d+u_{(d-1)y}z^{d-1}\\
+\cdots+u_{1y}z+u_{0y})-h_x(du_dz^{d-1}+\cdots+u_1)-v_1h_y=0.
\end{split}
\end{equation}
We always view that the polynomials are in $K[x,y,z]$ with coefficients in $K[x,y]$ when comparing the coefficients of the degree of $z$. Comparing the coefficients of $z^{d+1}$ and $z^d$ of equation \eqref{eq2.39}, we have the following equations:
\begin{equation}\label{eq2.40}
u_{dy}v_{1x}=0
\end{equation}
and
\begin{equation}\label{eq2.41}
v_{1x}u_{(d-1)y}+v_{0x}u_{dy}=0.
\end{equation}
It follows from equations \eqref{eq2.40} and \eqref{eq2.41} that $v_{1x}=v_{0x}=0$ or $v_{1x}=u_{dy}=0$ or $u_{dy}=u_{(d-1)y}=0$.\\

(a) If $v_{1x}=v_{0x}=0$, then equation \eqref{eq2.39} has the following
form:
\begin{equation}\label{eq2.42}
(u_{1x}z+u_{0x})(v_{1y}z+v_{0y})-h_x(du_dz^{d-1}+\cdots+u_1)-h_yv_1=0.
\end{equation}
If $d>3$, then $h_x=0$ by comparing the coefficient of $z^{d-1}$ of equation \eqref{eq2.42}.
Thus, it follows from equations \eqref{eq2.5} and \eqref{eq2.38} that
$h_y(u_{1x}z+u_{0x})=0$. Therefore, we have $h_y=0$ or
$u_{1x}=u_{0x}=0$. If $u_{1x}=u_{0x}=0$, then it follows from equation
\eqref{eq2.42} that $h_y=0$. Thus, we have $h=0$ because
$h(0,0)=0$. Therefore, $u,v,h$ are linearly dependent.\\
If $d=3$, then comparing the coefficients of $z^2,z$ and $z^0$ of equation
\eqref{eq2.42}, we have the following equations:
\begin{equation}\label{eq2.43}
u_{1x}v_{1y}-3u_3h_x = 0,
\end{equation}
\begin{equation}\label{eq2.44}
  u_{1x}v_{0y}-v_{1y}u_{0x}-2u_2h_x = 0
\end{equation}
and
\begin{equation}\label{eq2.45}
u_{0x}v_{0y}-u_1h_x-v_1h_y = 0.
\end{equation}
It follows from equations \eqref{eq2.5} and \eqref{eq2.38} that
\begin{equation}\label{eq2.46}
\begin{split}
-h_x(v_{1y}z+v_{0y})(3u_3z^2+2u_2z+u_1)+[h_x(u_{3y}z^3+u_{2y}z^2+u_{1y}z+u_{0y})\\
-h_y(u_{1x}z+u_{0x})]v_1=0.
\end{split}
\end{equation}
Comparing the coefficients of $z^3$ of the above equation, we have
$h_x(3v_{1y}u_3-u_{3y}v_1)=0$. Thus, we have $h_x=0$ or
$3u_3v_{1y}=v_1u_{3y}$.

(a1) If $h_x=0$, then it follows from equation \eqref{eq2.46} that
$h_y(u_{1x}z+u_{0x})=0$. That is, $h_y=0$ or $u_{1x}=u_{0x}=0$. If
$u_{1x}=u_{0x}=0$, then it follows from equation \eqref{eq2.42} that
$h_y=0$. Thus, we have $h=0$ because $h(0,0)=0$. Therefore, $u,v,h$
are linearly dependent.

(a2) If $3u_3v_{1y}=v_1u_{3y}$, then
\begin{equation}\label{eq2.47}
\frac{u_{3y}}{u_3}=3\frac{v_{1y}}{v_1}.
\end{equation}
If $v_{1y}=0$, then $u_{3y}=0$. It follows from equation \eqref{eq2.43} that
$h_x=0$. Thus, it follows from the arguments of (a1) that
$u,v,h$ are linearly dependent. We can assume that $u_{3y}v_{1y}\neq
0$. Then we have that $u_3=e^{\bar{d}(x)}v_1^3$ by integrating the
two sides of equation \eqref{eq2.47} with respect to $y$, where $\bar{d}(x)$
is a function of $x$. Since $u_3,v_1\in  K[x,y]$, we have
$e^{\bar{d}(x)}\in K(x)$.
That is,
\begin{equation}\label{eq2.48}
u_3=d(x)v_1^3
\end{equation}
with $d(x)\in  K(x)$, $d(x)\neq 0$. Let $d(x)=\frac{d_1(x)}{d_2(x)}$ with
$d_1(x),d_2(x)\in  K[x]$ and $d_1(x)\cdot d_2(x)\neq 0$.
Substituting equations \eqref{eq2.7} and \eqref{eq2.48} into equation \eqref{eq2.43}, we
have
\begin{equation}\label{eq2.49}
-3d_1(x)v_1^3h_x=d_2(x)v_{1y}^2
\end{equation}
Then we have $v_{1y}=0$ by comparing the degree of $y$ of equation \eqref{eq2.49}. It follows from equation \eqref{eq2.49} that $h_x=0$. Thus, it reduces to (a1).\\

(b) If $v_{1x}=u_{dy}=0$, then it follows from equation \eqref{eq2.38} that
$u_d\in K^*$.
Thus, we have the following equation:
\begin{equation}\label{eq2.50}
-v_{0x}u_{iy}-(i+1)u_{i+1}h_x=0
\end{equation}
by comparing the coefficients of $z^i$ of equation \eqref{eq2.39} for
$i=d-1,d-2,\ldots,3$. Comparing the coefficients of $z^2,z$ and
$z^0$ of equation \eqref{eq2.39}, we have the following equations:
\begin{equation}\label{eq2.51}
u_{1x}v_{1y}-v_{0x}u_{2y}-3u_3h_x = 0,
\end{equation}
\begin{equation}\label{eq2.52}
u_{1x}v_{0y}+v_{1y}u_{0x}-v_{0x}u_{1y}-2u_2h_x = 0
\end{equation}
and
\begin{equation}\label{eq2.53}
u_{0x}v_{0y}-v_{0x}u_{0y}-u_1h_x-v_1h_y = 0.
\end{equation}
It follows from equations \eqref{eq2.5} and \eqref{eq2.38} that
\begin{equation}\label{eq2.54}
\begin{split}
[v_{0x}h_y-h_x(v_{1y}z+v_{0y})](du_dz^{d-1}+(d-1)u_{d-1}z^{d-2}+\cdots+u_1)\\
+[h_x(u_{(d-1)y}z^{d-1}+\cdots+u_{1y}z+u_{0y})-h_y(u_{1x}z+u_{0x})]v_1=0.
\end{split}
\end{equation}
Then we have $h_xv_{1y}=0$ by comparing the coefficients of $z^d$ of equation
\eqref{eq2.54}. That is, $h_x=0$ or $v_{1y}=0$.

(b1) If $h_x=0$, then equation \eqref{eq2.54} has the following form:
\begin{equation}\label{eq2.55}
v_{0x}h_y(du_dz^{d-1}+(d-1)u_{d-1}z^{d-2}+\cdots+u_1)-h_y(u_{1x}z+u_{0x})v_1=0.
\end{equation}
Comparing the coefficients of $z^{d-1}$ of equation \eqref{eq2.55}, we have
that $v_{0x}h_y=0$. That is, $v_{0x}=0$ or $h_y=0$. If $v_{0x}=0$,
then it reduces to (a). If $h_y=0$, then $h=0$ because
$h(0,0)=0$. Thus, $u,v,h$ are linearly dependent.

(b2) If $v_{1y}=0$, then comparing the coefficients of $z^{d-1}$ and
$z^0$ of equation \eqref{eq2.54}, we have
\begin{equation}\label{eq2.56}
(v_{0x}h_y-h_xv_{0y})du_d+h_xu_{(d-1)y}v_1=0
\end{equation}
and
\begin{equation}\label{eq2.57}
(v_{0x}h_y-h_xv_{0y})u_1+(h_xu_{0y}-h_yu_{0x})v_1=0.
\end{equation}
It follows from equation \eqref{eq2.50}$(i=d-1)$ for $d>3$ and from equation \eqref{eq2.51} for $d=3$ that
\begin{equation}\label{eq2.58}
h_x=-\frac{1}{du_d}v_{0x}u_{(d-1)y}.
\end{equation}
Substituting equation \eqref{eq2.58} into equation \eqref{eq2.56}, we have the following equation:
$$v_{0x}[h_y-\frac{v_1}{d^2u_d^2}u_{(d-1)y}^2+\frac{1}{du_d}u_{(d-1)y}v_{0y}]=0$$
for $d\geq 3$. Thus, we have $v_{0x}=0$ or
\begin{equation}\label{eq2.59}
h_y=\frac{v_1}{d^2u_d^2}u_{(d-1)y}^2-\frac{1}{du_d}u_{(d-1)y}v_{0y}.
\end{equation}
If $v_{0x}=0$, then it reduces to (a). Otherwise,
substituting equations \eqref{eq2.58}, \eqref{eq2.59} into equation \eqref{eq2.53},we have
that
\begin{equation}\label{eq2.60}
u_{0x}v_{0y}-v_{0x}u_{0y}=-\frac{u_1}{du_d}v_{0x}u_{(d-1)y}+\frac{v_1^2}{d^2u_d^2}u_{(d-1)y}^2-\frac{v_1}{du_d}u_{(d-1)y}v_{0y}.
\end{equation}
Substituting equations \eqref{eq2.58}, \eqref{eq2.59} into equation \eqref{eq2.57}, we have that
\begin{equation}\label{eq2.61}
\frac{u_1v_1}{d^2u_d^2}u_{(d-1)y}^2v_{0x}-\frac{v_1^2}{d^2u_d^2}u_{(d-1)y}^2u_{0x}+\frac{v_1}{du_d}u_{(d-1)y}(u_{0x}v_{0y}-v_{0x}u_{0y})=0.
\end{equation}
Then we have $\frac{v_1^3}{d^3u_d^3}u_{(d-1)y}^3=0$ by substituting equations
\eqref{eq2.7}, \eqref{eq2.60} into equation \eqref{eq2.61}. That is,
$u_{(d-1)y}=0$. It follows from equation \eqref{eq2.50}$(i=d-1)$ that
$h_x=0$. Then it reduces to (b1).\\

(c) If $u_{dy}=u_{(d-1)y}=0$, then it follows from equations \eqref{eq2.4}
and \eqref{eq2.6} that
\begin{equation}\label{eq2.62}
\begin{split}
(u_{1x}z+u_{0x})(v_{1y}z+v_{0y})-(v_{1x}z+v_{0x})(u_{(d-2)y}z^{d-2}+\cdots\\
+u_{1y}z+u_{0y})-h_x(du_dz^{d-1}+\cdots+u_1)-h_yv_1=0.
\end{split}
\end{equation}
Comparing the coefficients of $z^j$ of equation \eqref{eq2.62} for
$j=d-2,\ldots,3$, we have the following equations:
\begin{equation}\label{eq2.63}
-v_{1x}u_{(j-1)y}-v_{0x}u_{jy}-(j+1)u_{j+1}h_x=0.
\end{equation}
Comparing the coefficients of $z^{d-1}, z^2, z$ and $z^0$ of equation \eqref{eq2.62}, we
have the following equations:
\begin{equation}\label{eq2.64}
  -v_{1x}u_{(d-2)y}-du_dh_x=0,
\end{equation}
\begin{equation}\label{eq2.65}
  u_{1x}v_{1y}-v_{1x}u_{1y}-v_{0x}u_{2y}-3u_3h_x = 0,
\end{equation}
\begin{equation}\label{eq2.66}
  u_{1x}v_{0y}+v_{1y}u_{0x}-v_{1x}u_{0y}-v_{0x}u_{1y}-2u_2h_x = 0
\end{equation}
and
\begin{equation}\label{eq2.67}
  u_{0x}v_{0y}-v_{0x}u_{0y}-u_1h_x-v_1h_y = 0.
\end{equation}
If $d=3$, then equations \eqref{eq2.63} and \eqref{eq2.64} are not available. If $d=4$, then equation \eqref{eq2.63} is not available.
It follows from equations \eqref{eq2.5} and \eqref{eq2.38} that
\begin{equation}\label{eq2.68}
\begin{split}
[(v_{1x}z+v_{0x})h_y-h_x(v_{1y}z+v_{0y})](du_dz^{d-1}+(d-1)u_{d-1}z^{d-2}+\cdots\\
+u_1)+[h_x(u_{(d-2)y}z^{d-2}+\cdots+u_{1y}z+u_{0y})-h_y(u_{1x}z+u_{0x})]v_1=0.
\end{split}
\end{equation}
Comparing the coefficients of $z^d$ and $z^{d-1}$ of equation \eqref{eq2.68}, we have the following equations:
\begin{equation}
\nonumber
  \left\{ \begin{aligned}
  du_d(v_{1x}h_y-h_xv_{1y})=0,~~~~~~~~~~~~~~~~~~~~~~~~~~~~~~~~~~~~~~\\
  (d-1)u_{d-1}(v_{1x}h_y-h_xv_{1y})+du_d(v_{0x}h_y-h_xv_{0y}) = 0. \\
                          \end{aligned} \right.
\end{equation}
That is,
\begin{eqnarray}\label{eq2.69}
 v_{1x}h_y-h_xv_{1y}=0,~~~
 v_{0x}h_y-h_xv_{0y}=0.
\end{eqnarray}
 Then equation \eqref{eq2.68} has the following form:
 \begin{equation}\label{eq2.70}
h_x(u_{(d-2)y}z^{d-2}+\cdots+u_{1y}z+u_{0y})-h_y(u_{1x}z+u_{0x})=0.
\end{equation}
Then we have $h_xu_{ky}=0$ by comparing the coefficients of $z^k$ of equation \eqref{eq2.70} for $2\leq k\leq d-2$.

If $d\geq 4$, then we have $h_x=0$ or $u_{(d-2)y}=\cdots=u_{2y}=0$.
If $u_{(d-2)y}=\cdots=u_{2y}=0$, then it follows from equation \eqref{eq2.64} that $h_x=0$. If $h_x=0$, then it follows from equation \eqref{eq2.70} that $h_y=0$ or $u_{1x}=u_{0x}=0$. If $u_{1x}=u_{0x}=0$, then
it follows from equation \eqref{eq2.64} that $v_{1x}=0$ or
$u_{(d-2)y}=0$. If $v_{1x}=0$, then it reduces to (b). If $u_{(d-2)y}=0$, then it follows from equation
\eqref{eq2.63} that $u_{(d-3)y}=\cdots=u_{2y}=0$. It follows from equation \eqref{eq2.65}
that $u_{1y}=0$. Then we have $u_{0y}=0$ by substituting the above equations into equation \eqref{eq2.66}. It follows from equation \eqref{eq2.67} that $h_y=0$.
Thus, we have $h=0$ because $h(0,0)=0$. Therefore, $u,v,h$ are linearly dependent.

If $d=3$, then equation \eqref{eq2.70} has the following form:
$$h_x(u_{1y}z+u_{0y})-h_y(u_{1x}z+u_{0x})=0.$$
That is,
\begin{eqnarray}\label{eq2.71}
  u_{1y}h_x-h_yu_{1x}=0,~~~
  u_{0y}h_x-h_yu_{0x}=0.
\end{eqnarray}
If $h=0$, then the conclusion follows. Suppose that $h\neq 0$ in the following arguments. It follows from equations \eqref{eq2.69} and \eqref{eq2.71} that
\begin{eqnarray}\label{eq2.72}
v_{1y}:v_{1x}=h_y:h_x=v_{0y}:v_{0x},~~~u_{1y}:u_{1x}=h_y:h_x=u_{0y}:u_{0x}.
\end{eqnarray}
Substituting \eqref{eq2.72} into equations \eqref{eq2.65}, \eqref{eq2.66}, \eqref{eq2.67} respectively, we have the following equations:
\begin{equation}\label{eq2.73}
v_{0x}u_{2y}+3u_3h_x=0,
\end{equation}
\begin{equation}\label{eq2.74}
2u_2h_x=0
\end{equation}
and
\begin{equation}\label{eq2.75}
u_1h_x+v_1h_y=0.
\end{equation}
It follows from equation \eqref{eq2.74} that $u_2=0$ or $h_x=0$. If $u_2=0$, then it follows from equation \eqref{eq2.73} that $h_x=0$. If $h_x=0$, then it follows from equation \eqref{eq2.75}
that $h_y=0$. Thus, we have $h=0$ because $h(0,0)=0$. This contradicts with our assumption.
\end{proof}

\begin{prop}
Let $H=(u(x,y,z),v(x,y,z),h(x,y))$ be a polynomial map with
$\deg_zv(x,y,z)=1$. Assume that $H(0)=0$ and $\deg_zu\geq 2$. If
$JH$ is nilpotent, then there exists $T\in \operatorname{GL}_3(K)$ such that
$$T^{-1}HT=(a_2(z)h(a_1(z)x+a_2(z)y)+c_1(z),-a_1(z)h(a_1(z)x+a_2(z)y)+c_2(z),0)$$
for some $a_i(z), c_i(z)\in K[z]$ and $f(t)\in K[z][t]$.
\end{prop}
\begin{proof}
It follow from Theorem \ref{thm2.3} that $u$, $v$, $h$ are linearly dependent over $K$. Then the conclusion follows from Theorem 7.2.25 in \cite{3} or Corollary 1.1 in \cite{18}.
\end{proof}

Next we only need to consider the polynomial maps $H$ and the components of $H$ are linearly independent over $K$.

\begin{lem}\label{lem2.4}
Let $H=(u,v,h)$ be a polynomial map in $K[x,y,z]$ with nilpotent Jacobian matrix. Let $d$ be the degree of $(u,v,h)$ with respect to $z$, and $(u_d,v_d,h_d)$ be the coefficient of $z^d$ of $(u,v,h)$. Then we can transform linearly to obtain $v_d\in K$, $u_d\in K[x_2]$, and the degree with respect to $z$ of $h$ unchanged.
\end{lem}
\begin{proof}
Taking coefficients of $z^d$ and $z^{2d}$ of the trace condition and the $2\times 2$ minors condition respectively, we obtain that $J_{x_1,x_2}(u_d,v_d)$ is nilpotent. By way of a linear transformation, we obtain that $J_{x_1,x_2}(u_d,v_d)$ is upper triangular. This yields the claims.
\end{proof}

\begin{thm}\label{thm2.5}
Let $H=(u(x,y,z),v(x,y,z),h(x,y))$ be a polynomial map with
$\deg_zv(x,y,z)\leq 1$. Assume that $H(0)=0$ and the components of
$H$ are linearly independent over $K$. If $JH$ is nilpotent,
then there exists $T\in \operatorname{GL}_3(K)$ such that
$T^{-1}HT$ has the form of Theorem 2.4 for $n=3$ in \cite{14}.
\end{thm}
\begin{proof}
Since $u,v,h$ are linearly independent, it follows from Theorem \ref{thm2.3} that $\deg_z u\leq 1$. Then it follows from Lemma \ref{lem2.4} that there exists $T_1\in \operatorname{GL}_3(K)$ such that $T_1^{-1}HT_1=(u_1z+u_0,v_1z+v_0,h(x,y))$ with $u_1\in K[x_2]$, $v_1\in K$, $u_0,~v_0\in K[x,y]$. Taking the coefficients of $z$ of the $2\times 2$ minors condition and $3\times 3$ minor condition of $J(T_1^{-1}HT_1)$ respectively, we obtain that $v_{0x}u_{1y}=0$ and $v_1h_xu_{1y}=0$. Thus, we have $u_1\in K$ or $v_1=0$
or $v_{0x}=h_x=0$. As for the two former two cases, there exists $T_2\in \operatorname{GL}_3(K)$ such that $T_2^{-1}T_1^{-1}HT_1T_2=(\bar{u}(x,y,z),\bar{v}(x,y),h(x,y))$. Thus, the conclusion follows from Theorem 2.4 in \cite{14}. If $v_{0x}=h_x=0$, then the determinant of $J(T_1^{-1}HT_1)$ is $v_1u_{0x}h_y$, which is 0. Thus, we have $u_{0x}=0$ or $h_y=0$. If $u_{0x}=0$, then we have $h_y=0$ by considering the $2\times 2$ minors condition of $J(T_1^{-1}HT_1)$. If $h_y=0$, then $h=0$ because $H(0)=0$. Thus, $u,v,h$ are linearly dependent over $K$. This contradicts with the condition that the components of $H$ are linearly independent over $K$.
\end{proof}

\section{A generalization of the form of $H$}

In the section, we first prove in Lemma \ref{lem3.2} that $\deg H_1^{(d)}\leq 1$, where $H_1^{(d)}$ is the leading homogeneous part with respect to $x_3,\ldots,x_n$ of $H_1$, where $H=(H_1(x_1,x_2,\ldots,x_n),b_3x_3+\cdots+b_nx_n+H_2^{(0)}(x_1,x_2),H_3(x_1,x_2),\ldots,H_n(x_1,x_2))$, $JH$ is nilpotent and the components of $H$ are linearly independent.
 Then we classify in Theorem \ref{thm3.3} all polynomial maps of the form $$H=(H_1(x_1,x_2,\ldots,x_n),b_3x_3+\cdots+b_nx_n+H_2^{(0)}(x_2),H_3(x_1,x_2),\ldots,H_n(x_1,x_2)),$$ where $JH$ is nilpotent and the components of $H$ are linearly independent.

\begin{lem}\label{lem3.1}
  Let $H$ be a polynomial map over $K$ of the form
  $$(H_1(x_1,x_2,\ldots,x_n),H_2(x_1,x_2,\ldots,x_n),H_3(x_1,x_2),\ldots,H_n(x_1,x_2))$$
where $H_2=b_3x_3+\cdots+b_nx_n+H_2^{(0)}(x_1,x_2)$, $b_3,\ldots,b_n\in K$, $H_2^{(0)}\in K[x_1,x_2]$. Write $h_2=b_3H_3+\cdots+b_nH_n$. If $JH$ is nilpotent, then
\begin{eqnarray}
H_{1x_1}+H_{2x_2}=0,\label{eq3.1}\\
(H_{2x_2})^2+H_{1x_2}H_{2x_1}+H_{1x_3}H_{3x_1}+\cdots+H_{1x_n}H_{nx_1}+h_{2x_2}=0,\label{eq3.2}\\
\begin{split}
H_{1x_3}(H_{2x_1}H_{3x_2}-H_{2x_2}H_{3x_1})+H_{1x_4}(H_{2x_1}H_{4x_2}-H_{2x_2}H_{4x_1})+\cdots\\
+H_{1x_n}(H_{2x_1}H_{nx_2}-H_{2x_2}H_{nx_1})-(H_{1x_1}h_{2x_2}-H_{1x_2}h_{2x_1})=0,\label{eq3.3}\\
\end{split}
\end{eqnarray}
\begin{equation}
\begin{split}
H_{1x_3}(H_{3x_1}h_{2x_2}-H_{3x_2}h_{2x_1})+H_{1x_4}(H_{4x_1}h_{2x_2}-H_{4x_2}h_{2x_1}+\cdots\\
+H_{1x_n}(H_{nx_1}h_{2x_2}-H_{nx_2}h_{2x_1})=0.\label{eq3.4}\\
\end{split}
\end{equation}
\end{lem}
\begin{proof}
Equation \eqref{eq3.1} follows from the fact that the trace of $JH$ is zero. Since the sum of the principal minor determinants of size 2 of $JH$ is zero as well, we deduce that
$$-H_{1x_1}H_{2x_2}+H_{1x_2}H_{2x_1}+H_{1x_3}H_{3x_1}+\cdots+H_{1x_n}H_{nx_1}+h_{2x_2}=0.$$
Adding equation \eqref{eq3.1}$H_{2x_2}$ times to it yields equation \eqref{eq3.2}. Since the sum of the principal minor determinants of size 3 of $JH$ is zero as well, we deduce that
\begin{equation}
\nonumber
\begin{split}
(H_{3x_1}H_{4x_2}-H_{4x_1}H_{3x_2})(b_4H_{1x_3}-b_3H_{1x_4})+\cdots
+(H_{3x_1}H_{nx_2}-\\H_{nx_1}H_{3x_2})(b_nH_{1x_3}-b_3H_{1x_n})+(H_{4x_1}H_{5x_2}-H_{5x_1}H_{4x_2})(b_5H_{1x_4}-b_4H_{1x_5})\\
+\cdots+(H_{(n-1)x_1}H_{nx_2}-H_{nx_1}H_{(n-1)x_2})(b_nH_{1x_{n-1}}-b_{n-1}H_{1x_n})=0.\\
\end{split}
\end{equation}
We view the above equation in another way; put the terms which contain $H_{1x_i}$ together in the above equation for $3\leq i\leq n$, so we have equation \eqref{eq3.4}.
\end{proof}

\begin{lem}\label{lem3.2}
Let $H$ be a polynomial map over $K$ of the form
$$(H_1(x_1,x_2,\ldots,x_n),H_2(x_1,x_2,\ldots,x_n),H_3(x_1,x_2),\ldots,H_n(x_1,x_2)),$$
where $H(0)=0$ and $H_2(x_1,\ldots,x_n)=b_3x_3+\cdots+b_nx_n+H_2^{(0)}(x_1,x_2)$, $b_3,\ldots,b_n\in K$, $H_2^{(0)}\in K[x_1,x_2]$. If $JH$ is nilpotent and the components of $H$ are linearly independent over $K$, then $\deg H_1^{(d)}\leq 1$, where $H_1^{(d)}$ is the leading homogeneous part with respect to $x_3,\ldots,x_n$ of $H_1$. Moreover, If $\deg H_1^{(d)}= 1$, then $H_1^{(d)}\in K[x_3,x_4,\ldots,x_n]$.
\end{lem}
\begin{proof}
If $b_3=\cdots=b_n=0$, then the conclusion follows from Lemma 2.3 in \cite{14}. Assume that at least one of $b_3,\ldots,b_n$ is non-zero in the following arguments. Write $H_1=H_1^{(d)}+H_1^{(d-1)}+\cdots+H_1^{(1)}+H_1^{(0)}$, where $H_1^{(i)}$ is the homogeneous part of degree $i$ with respect to $x_3,\ldots,x_n$ of $H_1$. Comparing the degree in equation \eqref{eq3.1} of the monomials with respect to $x_3,\ldots,x_n$ of degree $i$ for $0\leq i\leq d$, we have the following equations:
\begin{eqnarray}\label{eq3.5}
(H_1^{(d)})_{x_1}=\cdots=(H_1^{(1)})_{x_1}=0,~~~
(H_1^{(0)})_{x_1}+(H_2^{(0)})_{x_2}=0.
\end{eqnarray}

(a) If $H_{2x_1}=0$, then $(H_1^{(d)})_{x_2}\cdot h_{2x_1}=0$ by focusing on the leading homogeneous part with respect to $x_3, x_4, \ldots, x_n$ of equation \eqref{eq3.3}. Thus, we have $(H_1^{(d)})_{x_2}=0$ or $h_{2x_1}=0$.

If $h_{2x_1}=0$, then equation \eqref{eq3.4} has the following form:
$$h_{2x_2}[H_{1x_3}H_{3x_1}+H_{1x_4}H_{4x_1}+\cdots+H_{1x_n}H_{nx_1}]=0.$$
That means $h_{2x_2}=0$ or $H_{1x_3}H_{3x_1}+H_{1x_4}H_{4x_1}+\cdots+H_{1x_n}H_{nx_1}=0$. If $h_{2x_2}=0$, then $h=0$ because $H(0)=0$. Thus, $H_3,\ldots,H_n$ are linearly dependent. This contradicts the fact that the components of $H$ are linearly independent over $K$. Therefore, we have
\begin{equation}\label{eq3.6}
H_{1x_3}H_{3x_1}+H_{1x_4}H_{4x_1}+\cdots+H_{1x_n}H_{nx_1}=0.
\end{equation}
Substituting equation \eqref{eq3.6} into equations \eqref{eq3.2}, \eqref{eq3.3} respectively, we have the following equations
\begin{eqnarray}
(H_{2x_2})^2+h_{2x_2}=0,\label{eq3.7}\\
H_{1x_1}h_{2x_2}=0.\label{eq3.8}
\end{eqnarray}
Substituting equation \eqref{eq3.1} into equation \eqref{eq3.8}, we have the following equation:
\begin{equation}
H_{2x_2}h_{2x_2}=0.\label{eq3.9}
\end{equation}
It follows from equations \eqref{eq3.7} and \eqref{eq3.9} that $H_{2x_2}=h_{2x_2}=0$. Thus, we have $h=0$ because $H(0)=0$. Thus, $H_3,\ldots,H_n$ are linearly dependent. This contradicts the fact that the components of $H$ are linearly independent over $K$. Thus, we have that
$$(H_1^{(d)})_{x_2}=0.$$

(b) If $H_{2x_1}\neq 0$, then we have $(H_1^{(d)})_{x_2}=0$ by considering the leading homogeneous part with respect to $x_3, x_4, \ldots, x_n$ of equation \eqref{eq3.2}.\\

Now assume that $d>1$. Focus on the homogeneous part of degree $d-1$ with respect to $x_3, x_4, \ldots, x_n$ of equations \eqref{eq3.2}, \eqref{eq3.3} and \eqref{eq3.4} respectively. we deduce that
\begin{equation}\label{eq3.10}
(H_1^{(d-1)})_{x_2}H_{2x_1}+(H_1^{(d)})_{x_3}H_{3x_1}+\cdots+(H_1^{(d)})_{x_n}H_{nx_1}=0,
\end{equation}
\begin{equation}\label{eq3.11}
\begin{split}
H_{2x_1}((H_1^{(d)})_{x_3}H_{3x_2}+\cdots+(H_1^{(d)})_{x_n}H_{nx_2})-H_{2x_2}((H_1^{(d)})_{x_3}H_{3x_1}+\\
\cdots+(H_1^{(d)})_{x_n}H_{nx_1})+(H_1^{(d-1)})_{x_2}h_{2x_1}=0
\end{split}
\end{equation}
and
\begin{equation}\label{eq3.12}
\begin{split}
h_{2x_2}((H_1^{(d)})_{x_3}H_{3x_1}+\cdots+(H_1^{(d)})_{x_n}H_{nx_1})=\\
h_{2x_1}((H_1^{(d)})_{x_3}H_{3x_2}+\cdots+(H_1^{(d)})_{x_n}H_{nx_2}).
\end{split}
\end{equation}
As $(H_1^{(d-1)})_{x_1}=0$, we have $H_1^{(d-1)}\in K[x_2,x_3,\ldots,x_n]$. We view $H_1^{(d-1)}$ as a polynomial in $K[x_3,\ldots,x_n]$ with coefficients in $K[x_2]$, then we have the following equations:
\begin{equation}\label{eq3.13}
e_3H_{3x_1}+e_4H_{4x_1}+\cdots+e_nH_{nx_1}=-q(x_2)H_{2x_1},
\end{equation}
\begin{equation}\label{eq3.14}
H_{2x_1}(e_3H_{3x_2}+\cdots+e_nH_{nx_2})+q(x_2)h_{2x_1}=H_{2x_2}(e_3H_{3x_1}+\cdots+e_nH_{nx_1}),
\end{equation}
\begin{equation}\label{eq3.15}
h_{2x_2}(e_3H_{3x_1}+\cdots+e_nH_{nx_1})=h_{2x_1}(e_3H_{3x_2}+\cdots+e_nH_{nx_2})
\end{equation}
by comparing the coefficients of any monomials $x_3^{j_3}\cdots x_n^{j_n}$ with $j_3+\cdots+j_n=d-1$ of equations \eqref{eq3.10}, \eqref{eq3.11} and \eqref{eq3.12} respectively, where $q(x_2)\in K[x_2]$, $e_3,\ldots,e_n\in K$ and at least one of $e_3,\ldots,e_n$ is non-zero. Since $H_{2x_1}=(H_2^{(0)})_{x_1}$, we have that
\begin{equation}\label{eq3.16}
e_3H_3+e_4H_4+\cdots+e_nH_n=-q(x_2)H_2^{(0)}+g(x_2)
\end{equation}
by integrating the two sides of equation \eqref{eq3.13} with respect to $x_1$, where $g(x_2)\in K[x_2]$ and $g(0)=0$. Differentiating the two sides of equation \eqref{eq3.16} with respect to $x_2$, we have that
\begin{equation}\label{eq3.17}
e_3H_{3x_2}+e_4H_{4x_2}+\cdots+e_nH_{nx_2}=-q'(x_2)H_2^{(0)}-q(x_2)(H_2^{(0)})_{x_2}+g'(x_2).
\end{equation}
Since $H_{2x_1}=(H_2^{(0)})_{x_1}$, we have
\begin{equation}\label{eq3.18}
-q'(x_2)H_2^{(0)}(H_2^{(0)})_{x_1}+g'(x_2)(H_2^{(0)})_{x_1}+q(x_2)h_{2x_1}=0
\end{equation}
by substituting equation \eqref{eq3.17} into equation \eqref{eq3.14}. Thus, we have
\begin{equation}\label{eq3.19}
-\frac{1}{2}q'(x_2)(H_2^{(0)})^2+g'(x_2)H_2^{(0)}+q(x_2)h_2=\bar{g}(x_2)
\end{equation}
by integrating the two sides of equation \eqref{eq3.18} with respect to $x_1$, where $\bar{g}(x_2)\in K[x_2]$ and $\bar{g}(0)=0$. Substituting equations \eqref{eq3.13}, \eqref{eq3.17} into equation \eqref{eq3.15}, we have that
\begin{equation}\label{eq3.20}
q(x_2)(H_2^{(0)})_{x_1}h_{2x_2}=(q(x_2)(H_2^{(0)})_{x_2}+q'(x_2)H_2^{(0)}-g'(x_2))h_{2x_1}.
\end{equation}

If $q(x_2)=0$, then it follows from equation \eqref{eq3.13} that
\begin{equation}\label{eq3.21}
e_3H_{3x_1}+\cdots+e_nH_{nx_1}=0.
\end{equation}
Substituting equation \eqref{eq3.21} into equation \eqref{eq3.14}, we have that
$$H_{2x_1}(e_3H_{3x_2}+\cdots+e_nH_{nx_2})=0.$$
That is, $H_{2x_1}=0$ or $e_3H_{3x_2}+\cdots+e_nH_{nx_2}=0$.

If $e_3H_{3x_2}+\cdots+e_nH_{nx_2}=0$, then $e_3H_3+\cdots+e_nH_n=0$ because $H(0)=0$. Thus, $H_3,\ldots, H_n$ are linearly dependent. This contradicts the fact that the components of $H$ are linearly independent over $K$.

If $H_{2x_1}=0$, then we assume without loss of generality that
$$(H_1^{(d)})_{x_3}, (H_1^{(d)})_{x_4},\ldots,(H_1^{(d)})_{x_k}$$
are linearly independent over $K$, and $(H_1^{(d)})_{x_{k+1}}=(H_1^{(d)})_{x_{k+2}}=\cdots=(H_1^{(d)})_{x_n}=0$. It is easy to see that $k\geq 3$. Then $(H_1^{(d)})_{x_3}, (H_1^{(d)})_{x_4},\ldots,(H_1^{(d)})_{x_k}$ are linearly independent over $K(x_1,x_2)$ as well. So if we focus on the leading homogeneous part with respect to $x_3,x_4,\ldots,x_n$ of equation \eqref{eq3.4}, we infer that
$$H_{ix_1}h_{2x_2}-H_{ix_2}h_{2x_1}=0$$
for each $i\in \{3,4,\ldots,k\}$. Consequently, $H_i$ is algebraically dependent over $K$ on $h_2$ for each $i\in \{3,4,\ldots,k\}$, and there exists an $f\in K[x_1,x_2]$, such that $H_i,h_2\in K[f]$ for each $i\in \{3,4,\ldots,k\}$. So if we focus the leading homogeneous part with respect to $x_3,\ldots,x_n$ of equation \eqref{eq3.2}, $H_{3x_1},\ldots, H_{kx_1}$ are linearly dependent over $K(x_3,\ldots,x_n)$, and hence over $K$. Since the rank of the sub-matrix of rows $3,4,\ldots,k$ of $JH$ is 1, the rows of this sub-matrix are linearly dependent over $K$ along with the entries of the first column. This contradicts the fact that the components of $H$ are linearly independent over $K$. So we can assume that $q(x_2)\neq 0$ in the following arguments.\\

It follows from equation \eqref{eq3.18} that
\begin{equation}\label{eq3.22}
h_{2x_1}=(q(x_2))^{-1}(q'(x_2)H_2^{(0)}(H_2^{(0)})_{x_1}-g'(x_2)(H_2^{(0)})_{x_1}).
\end{equation}
Substituting equation \eqref{eq3.22} into equation \eqref{eq3.20}, we have that
\begin{equation}\label{eq3.23}
h_{2x_2}=(q(x_2))^{-2}(q(x_2)(H_2^{(0)})_{x_2}+q'(x_2)H_2^{(0)}-g'(x_2))(q'(x_2)H_2^{(0)}-g'(x_2)).
\end{equation}
Differentiating the two sides of equation \eqref{eq3.19} with respect to $x_2$, we have that
\begin{equation}\label{eq3.24}
\begin{split}
-\frac{1}{2}q''(x_2)(H_2^{(0)})^2-q'(x_2)H_2^{(0)}(H_2^{(0)})_{x_2}+g''(x_2)H_2^{(0)}\\ +g'(x_2)(H_2^{(0)})_{x_2}+q'(x_2)h_2+q(x_2)h_{2x_2}=\bar{g}'(x_2).
\end{split}
\end{equation}
We can get the equations about $h_2$ and $h_{2x_2}$  from equations \eqref{eq3.19} and \eqref{eq3.23} respectively, and then substituting them into equation \eqref{eq3.24}, we have the following equation:
\begin{equation}\label{eq3.25}
\begin{split}
(-\frac{1}{2}q(x_2)q''(x_2)+\frac{3}{2}(q'(x_2))^2)(H_2^{(0)})^2+(q(x_2)g''(x_2)-3g'(x_2)q'(x_2))H_2^{(0)}\\
=q(x_2)\bar{g}'(x_2)-q'(x_2)\bar{g}(x_2)-(g'(x_2))^2.
\end{split}
\end{equation}
If $(H_2^{(0)})_{x_1}=0$, then it follows from equation \eqref{eq3.22} that $h_{2x_1}=0$. Thus, we have that $H_3,\ldots,H_n$ are linearly dependent by following the arguments of Lemma \ref{lem3.2} (a). This contradicts the fact that the components of $H$ are linear independent over $K$. Thus, we have $(H_2^{(0)})_{x_1}\neq 0$.

Comparing the degree of $x_1$ of equation \eqref{eq3.25}, we have
\begin{equation}\label{eq3.26}
q(x_2)q''(x_2)=3(q'(x_2))^2
\end{equation}
and
\begin{equation}\label{eq3.27}
q(x_2)g''(x_2)=3g'(x_2)q'(x_2).
\end{equation}
Thus, we have $q'(x_2)=0$ by comparing the coefficients of the highest degree of $x_2$ of equation \eqref{eq3.26}. Therefore, it follows from equation \eqref{eq3.27} that $g''(x_2)=0$. Then equation \eqref{eq3.25} has the following form:
$$q(x_2)\bar{g}'(x_2)=(g'(x_2))^2.$$
Let $c:=q(x_2)\in {K}^*$. Since $H(0)=0$, we have $g(0)=\bar{g}(0)=0$. So we assume that $g(x_2)=\tilde{c}x_2$ for some $\tilde{c} \in K$. Then $\bar{g}(x_2)=\frac{\tilde{c}^2}{c}x_2$.
It follows from equation \eqref{eq3.19} that
\begin{equation}\label{eq3.28}
h_2=\frac{\tilde{c}}{c^2}(-c H_2^{(0)}+\tilde{c}x_2).
\end{equation}

If $\tilde{c}=0$, then it follows from equation \eqref{eq3.28} that $h_2=0$. Thus, $H_3,\ldots,H_n$ are linearly dependent. This contradicts the fact that the components of $H$ are linear independent over $K$.

If $\tilde{c}\neq 0$, then let $r=\frac{\tilde{c}}{c}\neq 0$, we have $b_3H_3+\cdots+b_nH_n+rH_2^{(0)}=r^2x_2$. Let
\[\bar{T}=\left(
  \begin{array}{ccccc}
    1 & 0 & 0 & \cdots & 0 \\
    0 & \frac{1}{r} & -\frac{b_3}{r}& \cdots & -\frac{b_n}{r}\\
    0 & 0 & 1 & \cdots & 0\\
    \vdots & \vdots &\vdots &\ddots &\vdots\\
    0 & 0 & 0 & \cdots & 1\\
  \end{array}
\right).\]
Then we have that $\bar{T}^{-1}H\bar{T}=(\bar{H}_1, r\cdot x_2, \bar{H}_3, \ldots, \bar{H}_n)$. Since $JH$ is nilpotent, we have that $J(\bar{T}^{-1}H\bar{T})$ is nilpotent. However, the element of the second row and the second column of the matrix $(J(\bar{T}^{-1}H\bar{T}))^m$ is $r^m$, which is not equal to zero. This contradicts the fact that the matrix $J(\bar{T}^{-1}H\bar{T})$ is nilpotent. Thus, $d\leq 1$. If $d=1$, then we have $H_1^{(d)}\in K[x_3\ldots,x_n]$ by following the former arguments of (a) and (b) Lemma \ref{lem3.2}.
\end{proof}

\begin{thm}\label{thm3.3}
Let $H$ be a polynomial map over $K$ of the form
$$(H_1(x_1,x_2,\ldots,x_n),H_2(x_2,\ldots,x_n), H_3(x_1,x_2),\ldots,H_n(x_1,x_2)),$$
where $H(0)=0$ and $H_2(x_2,\ldots,x_n)=b_3x_3+\cdots+b_nx_n+H_2^{(0)}(x_2)$, $b_3,\ldots,b_n\in K$, $H_2^{(0)}\in K[x_2]$. If $JH$ is nilpotent and the components of $H$ are linearly independent over $K$, then there exists a $T\in \operatorname{GL}_n(K)$ such that $T^{-1}\circ H\circ T$ be the form of Theorem 2.4 in \cite{14}.
\end{thm}
\begin{proof}
If $b_3=\cdots=b_n=0$, then the conclusion follows from Theorem 2.4 in \cite{14}. We can assume that at least one of $b_3,\ldots,b_n$ is non-zero in the following arguments.

It follows from Lemma \ref{lem3.2} that $\deg H_1^{(d)}\leq 1$, where $H_1^{(d)}$ is the leading homogeneous part with respect to $x_3,\ldots,x_n$ of $H_1$.

If $\deg H_1^{(d)}=0$, then let $\tilde{T}=P_n(1,2)$; $\tilde{T}^{-1}H\tilde{T}$ is of the form of Theorem 2.4 in \cite{14} and $J(\tilde{T}^{-1}H\tilde{T})$ is nilpotent. Thus, the conclusion follows from Theorem 2.4 in \cite{14}.

If  $\deg H_1^{(d)}= 1$, then let $H_1=a_3x_3+a_4x_4+\cdots+a_nx_n+H_1^{(0)}(x_1,x_2)$, where $H_1^{(0)}\in K[x_1,x_2]$.
Since $JH$ is nilpotent, it follows from Lemma \ref{lem3.1} that we have the following equations:
\begin{eqnarray}
(H_1^{(0)})_{x_1}+(H_2^{(0)})_{x_2}=0,\label{eq3.31}\\
(H_1^{(0)})_{x_1}(H_2^{(0)})_{x_2}-h_{1x_1}-h_{2x_2}=0,\label{eq3.32}\\
-(H_2^{(0)})_{x_2}h_{1x_1}-((H_1^{(0)})_{x_1}h_{2x_2}-(H_1^{(0)})_{x_2}h_{2x_1})=0,\label{eq3.33}\\
h_{2x_2}h_{1x_1}-h_{2x_1}h_{1x_2}=0\label{eq3.34}
\end{eqnarray}
where $h_1=\sum_{i=3}^na_iH_i$, $h_2=\sum_{i=3}^nb_iH_i$. Clearly, $h_1\cdot h_2\neq 0$. It follows from equation \eqref{eq3.34} that there exists $f\in K[x_1,x_2]$, such that $h_1, h_2\in K[f]$. We have the following equation:
\begin{equation}\label{eq3.35}
H_1^{(0)}=-(H_2^{(0)})'\cdot x_1+W(x_2)
\end{equation}
by integrating the two sides of equation \eqref{eq3.31} with respect to $x_1$, where $W(x_2)\in K[x_2]$. It follows from equation \eqref{eq3.32} that
\begin{equation}\label{eq3.36}
h_{1x_1}=-h_{2x_2}-[(H_2^{(0)})']^2.
\end{equation}
Replacing $h_{1x_1}$ with equation \eqref{eq3.36} in equation \eqref{eq3.33}, we have the following equation:
\begin{equation}\label{eq3.37}
2(H_2^{(0)})'h_{2x_2}+(H_1^{(0)})_{x_2}h_{2x_1}=-[(H_2^{(0)})']^3.
\end{equation}
Substituting equation \eqref{eq3.35} for $(H_1^{(0)})_{x_2}$ in equation \eqref{eq3.37}, we have the following equation:
\begin{equation}\label{eq3.38}
h_2'(f)[2(H_2^{(0)})'\cdot f_{x_2}-(H_2^{(0)})''\cdot x_1\cdot f_{x_1}+W'(x_2)f_{x_1}]=-[(H_2^{(0)})']^3.
\end{equation}
If $f_{x_1}=0$, then $h_{1x_1}=h_{2x_1}=0$. It follows from equation \eqref{eq3.33} that $(H_1^{(0)})_{x_1}\cdot h_{2x_2}\allowbreak
=0$. Thus, we have $(H_1^{(0)})_{x_1}=0$ or $h_{2x_2}=0$.

If $(H_1^{(0)})_{x_1}=0$, then it follows from equation \eqref{eq3.32} that $h_{2x_2}=0$.

If $h_{2x_2}=0$, then we have $h_2=0$ because $H(0)=0$. Thus, $H_3,\ldots,H_n$ are linearly dependent. This contradicts the fact that the components of $H$ are linearly independent over $K$.

If $f_{x_1}\neq 0$, then it follows from equation \eqref{eq3.38} that $h_2'(f)\in K$ or $(H_2^{(0)})'=0$ and $W'(x_2)=0$.

If $(H_2^{(0)})'=0$ and $W'(x_2)=0$, then we have $H_2^{(0)}=0$ and $H_1^{(0)}=W(x_2)=0$ because $H(0)=0$. Thus, it follows from equations \eqref{eq3.32} and \eqref{eq3.34} that $J(h_1,h_2)$ is nilpotent, which we deduce that there is $c\in K$ such that $h_2=ch_1$.

If $h_2'(f)\in K$, we have $h_2(f)=c_2f$ for some $c_2\in K$ because $H(0)=0$. It follows from equation \eqref{eq3.36} that
\begin{equation}\label{eq3.39}
h_1'(f)\cdot f_{x_1}=-c_2f_{x_2}-[(H_2^{(0)})']^2.
\end{equation}
Let $l=\deg_{x_1}f$. Then $l\geq 1$.

If $l\geq 2$, then $h_1'(f)\in K$ by comparing the degree of $x_1$ of equation \eqref{eq3.39}. Since $H(0)=0$, we have $h_1=c_1f$ for some $c_1\in K$.

If $l=1$, then let $f=\alpha_1(x_2)\cdot x_1+\alpha_0(x_2)$ with $\alpha_1, \alpha_0\in K[x_2]$ and $\alpha_1\neq 0$, we have $\deg_fh_1'\leq 1$ by comparing the degree of $x_1$ of equation
\eqref{eq3.39}. Let $h_1'=t_2f+c_1$ for some $c_1,t_2\in K$. We view that the polynomials are in $K[x_2][x_1]$ with coefficients in $K[x_2]$ when comparing the coefficients of $x_1^j$. Comparing the coefficients of $x_1$ of equation \eqref{eq3.39}, we have that
$$t_2\cdot \alpha_1^2=-c_2\alpha_1'.$$
Thus, we have that $\alpha_1'=0$ and $t_2=0$ by comparing the degree of $x_2$ of the above equation. Thus, we have
$$h_1=c_1f.$$
If $c_1=0$, then $h_1=0$. Thus, $H_3,\ldots,H_n$ are linearly dependent. This contradicts the fact that the components of $H$ are linearly independent over $K$.

If $c_1\neq 0$, then $h_2=\frac{c_2}{c_1}h_1$. Since the components of $H$ are linearly independent over $K$, we have $b_i=\frac{c_2}{c_1}a_i$ for all $3\leq i\leq n$. Let $\hat{T}=P_n(1(\frac{c_2}{c_1}),2)$. Then $\hat{T}^{-1}H\hat{T}$ is of the form of Theorem 2.4 in \cite{14}, and $J(\hat{T}^{-1}H\hat{T})$ is nilpotent. Thus, the conclusion follows.
\end{proof}

\begin{cor}
Let $F=x+H$, where $H$ be as in Theorem \ref{thm3.3}. If $JH$ is nilpotent and the components of $H$ are linearly independent over $K$, then $F$ is tame.
\end{cor}
\begin{proof}
The conclusion follows from Theorem \ref{thm3.3} and the arguments of section 3 in \cite{14}.
\end{proof}

\section{Some Remarks}

In order to classify all polynomial maps with nilpotent Jacobians of the form
$$(H_1(x_1,x_2,\ldots,x_n),H_2(x_1,x_2,\ldots,x_n),H_3(x_1,x_2),\ldots,H_n(x_1,x_2)),$$
where $H(0)=0$, $H_2(x_1,\ldots,x_n)=b_3x_3+\cdots+b_nx_n+H_2^{(0)}(x_1,x_2)$, and the components of $H$ are linearly independent over $K$, it suffices to classify all polynomial maps in dimension 4 of the form
$$\tilde{h}=(z+\tilde{h}_1(x,y),w+\tilde{h}_2(x,y),\tilde{h}_3(x,y),\tilde{h}_4(x,y)),$$
where $J\tilde{h}$ is nilpotent and $\tilde{h}_i\in K[x,y]$ for all $1\leq i\leq 4$, and the component of $\tilde{h}$ are linearly independent over $K$.

Jacobian nilpotency of $\tilde{h}$ is just the Keller condition on $x + t \tilde{h}$, so we have that
$$(x + t \tilde{h}_1 (x , y) + t z, y + t \tilde{h}_2 (x , y) + t w, z + t \tilde{h}_3 (x , y), w + t \tilde{h}_4 (x , y))$$
is a Keller map over $K[t]$. This is equivalent to that
$$(x + t \tilde{h}_1 (x , y) - t^2 \tilde{h}_3 (x , y), y + t \tilde{h}_2 (x , y) - t^2 \tilde{h}_4 (x , y))$$
is a Keller map over $K[t]$. By Moh's result \cite{Moh}, these maps are invertible over $K(t)$, hence over $K[t]$, if the degree is at most 100. In \cite{Wang}, it is claimed that there are errors in Moh's work, but these errors are repaired.

We can find the solution for $\tilde{h}$ if the components of $\tilde{h}$ are linearly dependent. More precisely, we have the following theorem.
\begin{thm}
Let $\tilde{h}=(z+\tilde{h}_1(x,y),w+\tilde{h}_2(x,y),\tilde{h}_3(x,y),\tilde{h}_4(x,y))$, where $\tilde{h}_i\in K[x,y]$ for $1\leq i\leq 4$. If $J\tilde{h}$ is nilpotent and the components of $\tilde{h}$ are linearly dependent, then there exists $T\in \operatorname{GL}_4(K)$ such that $T^{-1}\tilde{h}T=(0,w,0,0)+\tilde{H}$ and $\tilde{H}$ is the form of Theorem 3.1 in \cite{14}.
\end{thm}
\begin{proof}
Since the components of $\tilde{h}$ are linearly dependent, there exists $\lambda \in K$ such that $h_4=\lambda h_3$. Let $T_1\in P_4(3(\lambda),4)$. Then $T_1^{-1}\tilde{h}T_1=(z+\tilde{h}_1,w+\lambda z+\tilde{h}_2,\tilde{h}_3,0)$. There is $T_2\in P_4(1(\lambda),2)$ such that $T_2^{-1}(T_1^{-1}\tilde{h}T_1)T_2=(z+\hat{h}_1,w+\hat{h}_2-\lambda \hat{h}_1,\hat{h}_3,0)$, where $\hat{h}_i=\tilde{h}_i(x,y+\lambda x)$. Since $(T_1T_2)^{-1}\tilde{h}T_1T_2$ is nilpotent, we have $J\hat{H}$ is nilpotent, where $\hat{H}=(z+\hat{h}_1,\hat{h}_2-\lambda \hat{h}_1,\hat{h}_3,0)$. Thus, the conclusion follows from Theorem 3.1 in \cite{14}.
\end{proof}

{\bf{Acknowledgement}}: The author is very grateful to Michiel de Bondt who give some good suggestions, especially the proof of Theorem 2.5 and the setup of section 4.


\begin{thebibliography}{99}
\bibitem{2} H. Bass, E. Connell, D. Wright, \newblock {\em The Jacobian Conjecture: Reduction of Degree and Formal Expansion of the Inverse}, Bulletin of the American Mathematical Society, 7 (1982), 287-330.
\bibitem{10} M. de Bondt, \newblock {\em Quasi-translations and counterexamples to the homogeneous dependence problem}, Proceedings of the American Mathematical Society 134 (2006) 2849-2856.
\bibitem{HKM} M. de Bondt, \newblock {\em Homogeneous Keller maps}, Ph.D. thesis, Radboud University, Nijmegen, The Netherlands, July 2009.
\bibitem{B2} M. de Bondt, \newblock {\em Quadratic polynomial maps with Jacobian rank two}, Linear Algebra and its Applications, 565 (2019) 267-286.
\bibitem{8} M. de Bondt, A. van den Essen, \newblock {\em The Jacobian conjecture: linear triangularization for homogeneous polynomial maps in dimension three}, Report 0413, University of Nijmegen, The Netherlands, 2004.
\bibitem{12} M. de Bondt, D. Yan, \newblock {\em Triangularization properties of power linear maps and the structural conjecture}, Annales Polonici Mathematici, 112(3) (2014) 247-266.
\bibitem{CE} \'{A} Casta\~{n}eda, A. van den Essen, \newblock {\em A new class of Nilpotent Jacobians in any dimension}, arXiv:1804.00584.
\bibitem{18} M. Chamberland, A. van den Essen, \newblock {\em Nilpotent Jacobians in dimension three}, Journal of Pure and Applied Algebra, 205 (2006) 146-155.
\bibitem{9} A. van den Essen, \newblock {\em Nilpotent Jacobian matrices with independent rows}, Report 9603, University of Nijmegen, The Netherlands, 1996.
\bibitem{E} A. van den Essen, \newblock {\em A counterexample to Meisters' cubic-linear linearization Conjecture}, Indagationes Mathematicae, 9(3) (1998), 333-339.
\bibitem{3} A. van den Essen, \newblock {\em Polynomial Automophisms and the Jacobian Conjecture}, Vol.\@ 190 in Progress in Mathematics Birkhauser Basel, 2000.
\bibitem{7} E. Hubbers, \newblock {\em The Jacobian conjecture: cubic homogeneous maps in dimension four}, Master's Thesis, University of Nijmegen, The Netherlands, 1994.
\bibitem{1} O.H. Keller, \newblock {\em Ganze Cremona-transformationen}, Monatshefte F\"{u}r Mathematik Und Physik, 47(1) (1939) 229-306.
\bibitem{Moh} T.T.Moh, \newblock {\em On the Jacobian conjecture and the configuration of roots}, J. Reine Angew. Math. 340 (1983), 140-212.
\bibitem{SFGZ} R. dos Santos Freire Jr., G. Gorni, G. Zampieri, \newblock {\em Search for homogeneous polynomial invariants and a cubic-homogeneous mapping without quadratic invariants}, Universitatis Iagellonicae Acta Mathematica, 46 (2008), 7-13.
\bibitem{XS5} X. Sun, \newblock {\em Classification of Quadratic Homogeneous Automorphisms in Dimension Five}, Communications in Algebra, 42(7) (2014), 2821-2840.
\bibitem{Wang} Lih-Chung Wang, \newblock {\em On the Jacobian conjecture}, Taiwanese Journal of Mathematics, 9(2005), 421-431.
\bibitem{4} S.S.-S. Wang, \newblock {\em A Jacobian criterion for separability}, Journal of Algebra, 65 (1980), 453-494.
\bibitem{6} D. Wright, \newblock {\em The Jacobian conjecture: linear triangularization for cubics in dimension three}, Linear and Multilinear Algebra 34 (1993) 85-97.
\bibitem{5} A.V. Yagzhev, \newblock {\em On Keller's problem}, Siberian Mathematical Journal, 21 (1980), 747-754.
\bibitem{14} D. Yan, M. de Bondt, \newblock {\em The classification of some polynomial maps with nilpotent Jacobians}, Linear Algebra and its Applications, 565 (2019), 287-308.
\bibitem{13} D. Yan, G. Tang, \newblock {\em Polynomial maps with nilpotent Jacobians in dimension three}, Linear Algebra and its Applications, 489 (2016), 298-323.

\end{thebibliography}
\end{document}